# On the Algebra of Derivation Groups of Some Low-Dimensional Leibniz Algebras

L. A. Kurdachenko, M. M. Semko, I. Ya. Subbotin

Let L be an algebra over a field F with the binary operations + and [-,-]. Then L is called a *left Leibniz algebra* if it satisfies the left Leibniz identity

$$[[a,b], c] = [a, [b,c]] - [b,[a,c]]$$

for all a, b, c ∈ L. We will also use another form of this identity:

$$[a, [b, c]] = [[a, b], c]+[b, [a, c]].$$

Leibniz algebras first appeared in the paper of A. Bloh [**BA1965**], but the term "Leibniz algebra" appears in the book of J.L. Loday [**LJ1992**] and his article [**LJ1993**]. In [**LP1993**], J. Lodey and T. Pirashvili conducted an in-depth study on Leibniz algebras' properties. The theory of Leibniz algebras has developed very intensely in many different directions. Some of the results of this theory were presented in the book [**AOR2020**]. Note that Lie algebras present a partial case of Leibniz algebras. Conversely, if L is a Leibniz algebra in which [a, a] = 0 for every element a ∈ L, then it is a Lie algebra. Thus, Lie algebras can be characterized as anticommutative Leibniz algebras.

As in Lie algebras, the structure of Leibniz algebras is greatly influenced by their derivation algebras.

Denote by **End**$_F$(L) the set of all linear transformations of L. Then L is an associative algebra by the operations + and o. As usual, **End**$_F$(L) is a Lie algebra by the operations + and [,], where [f, g] = f o g – g o f for all f, g ∈ **End**$_F$(L).

A linear transformation f of a Leibniz algebra L is called a *derivation* if

$$f([a, b]) = [f(a), b] + [a, f(b)] \text{ for all } a, b \in L.$$

Let **Der**(L) be the subset of all derivations of L. It is possible to prove that **Der**(L) is a subalgebra of the Lie algebra **End**$_F$(L). **Der**(L) is called the *algebra of derivations* of the Leibniz algebra L.

The influence on the structure of the Leibniz algebra of their algebras of derivations can be observed in the following result: If A is an ideal of a Leibniz algebra, then the factor – algebra of L by the annihilator of A is isomorphic to some subalgebra of **Der**(A) [**KOP2016**, Proposition 3.2]. Therefore, finding out the structure of algebras of derivations of Leibniz algebras is one of the important steps in the process of studying the structure of Leibniz algebras. It is natural to start studying the algebra of derivations of a Leibniz algebra, the structure of which has been studied quite fully. The structure of the algebra of derivations of finite dimensional one–generator Leibniz algebras was described in paper [**KSEY2021, SSY2021**], and infinite-dimensional one-generator Leibniz algebras in paper [**KSUY2022**]. Naturally, the question about the algebras of derivations of Leibniz algebras of small dimensions arises. In contrast to the Lie algebra, the situation with the Leibniz algebras of dimension 3 is very diverse. The Leibniz algebras of dimension 3 are described, and their most detailed description can be found in [**KPS2022**]. The papers [**KSEY2023, KSEY2023A**] were the first to describe the algebras of derivations nilpotent Leibniz algebras of dimension 3. More concretely, they have described the algebras of derivations of nilpotency class 3 Leibniz algebras, nilpotency class 2 Leibniz algebras with a center that has dimension 2, and the algebra of derivations of the nilpotent Leibniz algebras of dimension 3 of nilpotency class 2, the center of which has dimension 1.

Let L be a nilpotent Leibniz algebra of nilpotency class 2 whose center of L has dimension 1. Of course, we suppose that L is not a Lie algebra. Then there is an element $a_1$ such that $[a_1, a_1] = a_3 \neq 0$. Since factor-algebra $L/\zeta(L)$ is abelian, $a_3 \in \zeta(L)$. It follows that $[a_1, a_3] = [a_3, a_1] = [a_3, a_3] = 0$. Then $\zeta(L) = Fa_3$. For every element $x \in L$, we have: $[a_1, x], [x, a_1] \in \zeta(L) \leq < a_1 > = Fa_1 \oplus Fa_3$. It follows that a subalgebra $< a_1 >$ is an ideal of L. Since $\dim_F(< a_1 >) = 2$, $< a_1 > \neq L$. Choose an element b such that $b \notin < a_1 >$. We have $[b, a_1] = \gamma a_3$ for some $\gamma \in F$. If $\gamma \neq 0$, then put $b_1 = \gamma^{-1}b - a_1$. Then $[b_1, a_1] = 0$. The choice of $b_1$ shows that $b_1 \notin < a_1 >$. It follows that subalgebra $\mathbf{Ann}^{left}_L(a_1)$ has dimension 2. The paper [KSeY2023] considered the case when $\mathbf{Ann}^{left}_L(a_1)$ is an abelian subalgebra. The next natural step is the case when $\mathbf{Ann}^{left}_L(a_1)$ is a non-abelian subalgebra. In this case, $[x, x] \neq 0$ for each element $x \in \mathbf{Ann}^{left}_L(a_1)$ such that $x \notin \zeta(L)$. It follows that subalgebra $\mathbf{Ann}^{left}_L(a_1)$ is a one-generator nilpotent algebra of dimension 2. Moreover, because $[L, L] = \zeta(L) \leq \mathbf{Ann}^{left}_L(a_1)$, the latter one is an ideal of L. Let b be an element that generates $\mathbf{Ann}^{left}_L(a_1)$. Since $\mathbf{Ann}^{left}_L(a_1)$ is non-abelian, $b \notin \zeta(L)$. We have $[b, a_1] = 0$ and $[a_1, b] = \gamma a_3$ for some scalar $\gamma \in F$.

If $\gamma = 0$, then the fact that $\mathbf{Ann}^{left}_L(a_1) = < b >$ shows that $\mathbf{Ann}^{left}_L(a_1) = \mathbf{Ann}_L(a_1)$, so that $[< a_1 >, < b >] = <0>$.

If $\gamma \neq 0$, then put $a_2 = \gamma^{-1}b$ and obtain $[a_1, a_2] = a_3$. Clearly, $< a_2 > = < b > = \mathbf{Ann}^{left}_L(a_1)$. Furthermore, $[a_2, a_2] = \lambda a_3$ where $0 \neq \lambda \in F$.

Thus, we come to the following two types of nilpotent Leibniz algebra:

$\mathbf{Lei}_4(3, F) = Fa_1 \oplus Fa_2 \oplus Fa_3$ where $[a_1, a_1] = a_3$, $[a_2, a_2] = \lambda a_3, 0 \neq \lambda \in F$, $[a_1, a_3] = [a_3, a_1] = [a_3, a_3] = [a_2, a_3] = [a_3, a_2] = 0$, $[a_1, a_2] = [a_2, a_1] = 0$

and

$\mathbf{Lei}_5(3, F) = Fa_1 \oplus Fa_2 \oplus Fa_3$ where $[a_1, a_1] = [a_1, a_2] = a_3$, $[a_2, a_2] = \lambda a_3, 0 \neq \lambda \in F$, $[a_1, a_3] = [a_3, a_1] = [a_3, a_3] = [a_2, a_3] = [a_3, a_2] = 0$, $[a_2, a_1] = 0$

In other words, $\mathbf{Lei}_4(3, F) = L$ is the sum of the ideals $A_1 = Fa_1 \oplus Fa_3$ and $A_2 = Fa_2 \oplus Fa_3$, and moreover, $A_1, A_2$ are nilpotent one-generator Leibniz algebras of dimension 2, $[A_1, A_2] = [A_2, A_1] = <0>$, $\mathbf{Leib}(L) = [L, L] = Fa_3$, $\zeta^{left}(L) = \zeta^{right}(L) = \zeta(L) = Fa_3$.

Leibniz algebra $L_7$ is a sum of the ideals $A_1 = Fa_1 \oplus Fa_3$ and $A_2 = Fa_2 \oplus Fa_3$. Moreover, $A_1$, $A_2$ are nilpotent cyclic Leibniz algebras of dimension 2, $[A_1, A_2] = Fa_3$, $[A_2, A_1] = <0>$, **Leib**$(L_7) = [L_7, L_7] = Fa_3$, $\zeta^{\text{left}}(L_7) = \zeta^{\text{right}}(L_7) = \zeta(L_7) = Fa_3$.

We dedicate the current article to the description of the algebras of derivations of nilpotent Leibniz algebras **Lei**$_4$(3, F) and **Lei**$_5$(3, F).

Let x be an arbitrary element of **Lei**$_4$(3, F), $x = \xi_1 a_1 + \xi_2 a_2 + \xi_3 a_3$. We have:

$$[x, x] = [\xi_1 a_1 + \xi_2 a_2 + \xi_3 a_3, \xi_1 a_1 + \xi_2 a_2 + \xi_3 a_3] =$$
$$\xi_1^2 [a_1, a_1] + \xi_2^2 [a_2, a_2] = (\xi_1^2 + \lambda \xi_2^2) a_3.$$

If we suppose that $[x, x] = 0$, then we are led to the Leibniz algebras whose algebra of derivations has been considered in the papers [**KSEY2022, KSEY2023A**]. Therefore, we will suppose that $[x, x] \neq 0$. If $\xi_1 = 0$ or $\xi_2 = 0$, then clearly, $[x, x] \neq 0$. Suppose that $\xi_1 \neq 0$ and $\xi_2 \neq 0$. Then we can see that polynomial $X^2 + \lambda$ has no root in field F.

We say that a field F is *2–closed* if equation $x^2 = a$ has a solution in F for every element $a \neq 0$.

We note that every locally finite (in particular, finite) field of characteristic 2 is 2–closed.

Thus, we can see that a Leibniz algebra of type **Lei**$_4$(3, F) over a 2–closed field F cannot exist.

A Leibniz algebra L is said to be *extraspecial* if $[L, L] = \zeta(L)$ is an ideal of dimension 1.

Thus, we can see that Leibniz algebra **Lei**$_4$(3, F) is extraspecial.

We begin with some general properties of the algebra of derivations of Leibniz algebras. Here, we show some basic elementary properties of derivations that have been proven in the paper [**KSUY2022**]. Let us recall some definitions.

Every Leibniz algebra L has a specific ideal. Denote by **Leib**(L) the subspace generated by the elements $[a, a]$, $a \in L$. It is possible to prove that **Leib**(L) is an ideal of L. The ideal **Leib**(L) is called the *Leibniz kernel* of algebra L. By its definition, factor–algebra L/**Leib**(L) is a Lie algebra. And conversely, if K is an ideal of L such that L/K is a Lie algebra, then K includes a Leibniz kernel.

Let L be a Leibniz algebra. Define the lower central series of L,
$L = \gamma_1(L) \geq \gamma_2(L) \geq \ldots \geq \gamma_\alpha(L) \geq \gamma_{\alpha+1}(L) \geq \ldots \gamma_\delta(L)$,

by the following rule: $\gamma_1(L) = L$, $\gamma_2(L) = [L, L]$, and recursively, $\gamma_{\alpha+1}(L) = [L, \gamma_\alpha(L)]$ for every ordinal $\alpha$, and $\gamma_\lambda(L) = \bigcap_{\mu<\lambda}\gamma_\mu(L)$ for every limit ordinal $\lambda$. The last term $\gamma_\delta(L) = \gamma_\infty(L)$ is called the *lower hypocenter* of L. We have: $\gamma_\delta(L) = [L, \gamma_\delta(L)]$.

If $\alpha = k$ is a positive integer, then $\gamma_k(L) = [L, [L, [L, \ldots, L] \ldots L]$ is the *left–normed commutator* of k copies of L.

As usual, we say that a Leibniz algebra L is called *nilpotent* if there exists a positive integer k such that $\gamma_k(L) = <0>$. More precisely, L is said to be *nilpotent of nilpotency class c* if $\gamma_{c+1}(L) = <0>$ but $\gamma_c(L) \neq <0>$.

The *left* (respectively *right*) *center* $\zeta^{\text{left}}(L)$ (respectively $\zeta^{\text{right}}(L)$) of a Leibniz algebra L is defined by the rule below:

$$\zeta^{\text{left}}(L) = \{ x \in L \mid [x, y] = 0 \text{ for each element } y \in L \}$$

(respectively

$$\zeta^{\text{right}}(L) = \{ x \in L \mid [y, x] = 0 \text{ for each element } y \in L \} ).$$

It is not hard to prove that the left center of L is an ideal, but such is not true for the right center. Moreover, **Leib**(L) $\leq \zeta^{\text{left}}(L)$, so that $L/\zeta^{\text{left}}(L)$ is a Lie algebra. The right center is a subalgebra of L; in general, the left and right centers are different; they may even have different dimensions (see [**KOP2016**]).

$\zeta(L)$ , the *center* of L, is defined by the rule below:

$$\zeta(L) = \{ x \in L \mid [x, y] = 0 = [y, x] \text{ for each element } y \in L \}.$$

The center is an ideal of L.

Starting with the center, we can construct the upper central series of Leibniz algebra L

$$<0> = \zeta_0(L) \leq \zeta_1(L) \leq \zeta_2(L) \leq \ldots \leq \zeta_\alpha(L) \leq \zeta_{\alpha+1}(L) \leq \ldots \zeta_\gamma(L) = \zeta_\infty(L),$$

by the following rule: $\zeta_1(L) = \zeta(L)$ is the center of L, and recursively, $\zeta_{\alpha+1}(L)/\zeta_\alpha(L) = \zeta(L/\zeta_\alpha(L))$ for every ordinal $\alpha$, and $\zeta_\lambda(L) = \cup_{\mu<\lambda}\zeta_\mu(L)$ for every limit ordinal $\lambda$. By definition, each term of this series is an ideal of L. The last term $\zeta_\infty(L)$ of this series is called the *upper hypercenter* of L. If $L = \zeta_\infty(L)$, then L is called a *hypercentral* Leibniz algebra.

**Lemma 1.** *Let L be a Leibniz algebra over a field F and f be a derivation of L. Then $f(\zeta^{\text{left}}(L) \leq \zeta^{\text{left}}(L)$, $f(\zeta^{\text{right}}(L)) \leq \zeta^{\text{right}}(L)$, and $f(\zeta(L)) \leq \zeta(L)$.*

**1.2. Corollary.** *Let L be a Leibniz algebra over a field F and f be a derivation of L. Then $f(\zeta_\alpha(L)) \leq \zeta_\alpha(L)$ for every ordinal $\alpha$.*

Denote by $\Xi$ the classic monomorphism of **End**(L) in $\mathbf{M}_3(F)$ (i.e., the mapping), assigning to each endomorphism its matrix with respect to the basis $\{ a_1, a_2, a_3 \}$.

**Theorem 1.** *Let D be an algebra of derivations of the Leibniz algebra $\mathbf{Lei}_4(3, F)$.*

*If **char**(F) = 2, then D is isomorphic to a Lie subalgebra of $\mathbf{M}_3(F)$ consisting of the following form's matrices:*

$$\begin{pmatrix} \alpha & 0 & 0 \\ 0 & \beta & 0 \end{pmatrix}$$

$\begin{pmatrix} & & \\ \gamma & \kappa & 0 \end{pmatrix}$, $\alpha$, $\beta$, $\gamma$, $\kappa \in F$.

*Furthermore, D is a direct sum of an abelian ideal **W** + **Z** and abelian subalgebra **U** + **V** where **Z** is the subspace of **D** generated by a derivation **z** defined by the rule $z(a_1) = a_3$, $z(a_2) = 0$, $z(a_3) = 0$; **W** is the subspace of **D** generated by a derivation **w** defined by the rule $w(a_1) = 0$, $w(a_2) = a_3$, $w(a_3) = 0$; **U** is the subspace of **D** generated by a derivation **u** defined by the rule $u(a_1) = a_1$, $u(a_2) = 0$, $u(a_3) = 0$; **V** is the subspace of **D** generated by a derivation **v** defined by the rule $v(a_1) = 0$, $v(a_2) = a_2$, $v(a_3) = 0$. Moreover, [**u**, **z**] = − **z**, [**u**, **w**] = − **w**, [**v**, **z**] = **0**, [**v**, **w**] = − **w**.*

*If **char**(F) $\neq$ 2, then D is isomorphic to a Lie subalgebra of **M**$_3$(F) consisting of the matrices of the following form:*

$\begin{pmatrix} 0 & 0 & 0 \\ 0 & 0 & 0 \\ \alpha & \beta & 0 \end{pmatrix}$, $\alpha$, $\beta \in F$.

*Furthermore, D is a direct sum of the subalgebras **W** and **Z** where **Z** is the subspace of **D** generated by a derivation **z** defined by the rule $z(a_1) = a_3$, $z(a_2) = 0$, $z(a_3) = 0$; **W** is the subspace of **D** generated by a derivation **w** defined by the rule $w(a_1) = 0$, $w(a_2) = a_3$, $w(a_3) = 0$; in particular, **D** is abelian.*

PROOF. Let L = **Lei**$_4$(3, F) and f $\in$ **Der**(L). By Lemma 1, $f(Fa_3) \leq Fa_3$, so that $f(a_1) = \alpha_1 a_1 + \alpha_2 a_2 + \alpha_3 a_3$, $f(a_2) = \beta_1 a_1 + \beta_2 a_2 + \beta_3 a_3$, $f(a_3) = \gamma a_3$. Then

$f(a_3) = f([a_1, a_1]) = [f(a_1), a_1] + [a_1, f(a_1)] =$
$[\alpha_1 a_1 + \alpha_2 a_2 + \alpha_3 a_3, a_1] + [a_1, \alpha_1 a_1 + \alpha_2 a_2 + \alpha_3 a_3] = \alpha_1[a_1, a_1] + \alpha_1[a_1, a_1] = 2\alpha_1 a_3$,
$f(a_3) = f(\lambda^{-1}[a_2, a_2]) = \lambda^{-1} f([a_2, a_2]) = \lambda^{-1}[f(a_2), a_2]] + \lambda^{-1}([a_2, f(a_2)] =$
$\lambda^{-1} [\beta_1 a_1 + \beta_2 a_2 + \beta_3 a_3, a_2] + \lambda^{-1} [a_2, \beta_1 a_1 + \beta_2 a_2 + \beta_3 a_3] =$
$\lambda^{-1} \beta_2[a_2, a_2] + \lambda^{-1} \beta_2[a_2, a_2] = 2\lambda^{-1} \beta_2\lambda a_3] = 2\beta_2 a_3$,
$0 = f([a_1, a_2]) = [f(a_1), a_2] + [a_1, f(a_2)] =$
$[\alpha_1 a_1 + \alpha_2 a_2 + \alpha_3 a_3, a_2] + [a_1, \beta_1 a_1 + \beta_2 a_2 + \beta_3 a_3] = \beta_1 a_3$.

Thus, we obtain the following equalities:

$2\alpha_1 = 2\beta_2 = \beta_1 = \gamma$.

In particular, if **char**(F) $\neq$ 2, then $\alpha_1 = \beta_2$.
If **char**(F) = 2, then $\Xi(f)$ is the following matrix:

$\begin{pmatrix} \alpha_1 & 0 & 0 \\ \alpha_2 & \beta_2 & 0 \\ \alpha_3 & \beta_3 & 0 \end{pmatrix}$, $\alpha_1$, $\alpha_2$, $\alpha_3$, $\beta_2$, $\beta_3 \in F$.

If **char**(F) $\neq$ 2, then $\Xi(f)$ is the following matrix:
$\begin{pmatrix} \alpha_1 & 2\alpha_1 & 0 \\ \alpha_2 & \alpha_1 & 0 \\ \alpha_3 & \beta_3 & 2\alpha_1 \end{pmatrix}$, $\alpha_1$, $\alpha_2$, $\alpha_3$, $\beta_1$, $\beta_2$, $\beta_3$, $\gamma \in F$.

Conversely, suppose that **char**(F) = 2 and let **f** be the linear transformation of L such that

$$\Xi(f) = \begin{pmatrix} \alpha_1 & 0 & 0 \\ \alpha_2 & \beta_2 & 0 \\ \alpha_3 & \beta_3 & 0 \end{pmatrix}, \ \alpha_1, \alpha_2, \alpha_3, \beta_2, \beta_3 \in F.$$

Let $x$, $y$ be arbitrary elements of $L$, $x = \xi_1 a_1 + \xi_2 a_2 + \xi_3 a_3$, $y = \eta_1 a_1 + \eta_2 a_2 + \eta_3 a_3$ where $\xi_1, \ \xi_2, \xi_3, \eta_1, \eta_2, \eta_3$ are arbitrary scalars. Then

$[x, y] = [\xi_1 a_1 + \xi_2 a_2 + \xi_3 a_3, \eta_1 a_1 + \eta_2 a_2 + \eta_3 a_3] = \xi_1 \eta_1 [a_1, a_1] + \xi_2 \eta_2 [a_2, a_2] = \xi_1 \eta_1 a_3 + \lambda \xi_2 \eta_2 a_3 = (\xi_1 \eta_1 + \lambda \xi_2 \eta_2) a_3;$

$\mathbf{f}(x) = \mathbf{f}(\xi_1 a_1 + \xi_2 a_2 + \xi_3 a_3) = \xi_1 \mathbf{f}(a_1) + \xi_2 \mathbf{f}(a_2) + \xi_3 \mathbf{f}(a_3) = \xi_1 (\alpha_1 a_1 + \alpha_2 a_2 + \alpha_3 a_3) + \xi_2(\beta_2 a_2 + \beta_3 a_3) = \xi_1 \alpha_1 a_1 + \xi_1 \alpha_2 a_2 + \xi_1 \alpha_3 a_3 + \xi_2 \beta_2 a_2 + \xi_2 \beta_3 a_3 = \xi_1 \alpha_1 a_1 + (\xi_1 \alpha_2 + \xi_2 \beta_2) a_2 + (\xi_1 \alpha_3 + \xi_2 \beta_3) a_3;$

$\mathbf{f}(y) = \eta_1 \alpha_1 a_1 + (\eta_1 \alpha_2 + \eta_2 \beta_2) a_2 + (\eta_1 \alpha_3 + \eta_2 \beta_3) a_3;$

$\mathbf{f}([x, y]) = \mathbf{f}((\xi_1 \eta_1 + \lambda \xi_2 \eta_2) a_3) = (\xi_1 \eta_1 + \lambda \xi_2 \eta_2) \mathbf{f}(a_3) = 0;$

$[\mathbf{f}(x), y] + [x, \mathbf{f}(y)] = [\xi_1 \alpha_1 a_1 + (\xi_1 \alpha_2 + \xi_2 \beta_2) a_2 + (\xi_1 \alpha_3 + \xi_2 \beta_3) a_3, \eta_1 a_1 + \eta_2 a_2 + \eta_3 a_3] + [\xi_1 a_1 + \xi_2 a_2 + \xi_3 a_3, \eta_1 \alpha_1 a_1 + (\eta_1 \alpha_2 + \eta_2 \beta_2) a_2 + (\eta_1 \alpha_3 + \eta_2 \beta_3) a_3] = \xi_1 \alpha_1 \eta_1 [a_1, a_1] + (\xi_1 \alpha_2 + \xi_2 \beta_2) \eta_2 [a_2, a_2] + \xi_1 \eta_1 \alpha_1 [a_1, a_1] + \xi_2 (\eta_1 \alpha_2 + \eta_2 \beta_2) [a_2, a_2] = \xi_1 \alpha_1 \eta_1 a_3 + (\xi_1 \alpha_2 + \xi_2 \beta_2) \eta_2 \lambda a_3 + \xi_1 \eta_1 \alpha_1 a_3 + \xi_2 (\eta_1 \alpha_2 + \eta_2 \beta_2) \lambda a_3 = (\xi_1 \alpha_1 \eta_1 + \lambda \xi_1 \alpha_2 \eta_2 + \lambda \xi_2 \beta_2 \eta_2 + \xi_1 \eta_1 \alpha_1 + \lambda \xi_2 \eta_1 \alpha_2 + \lambda \xi_2 \eta_2 \beta_2) a_3 = (2\xi_1 \alpha_1 \eta_1 + \lambda \alpha_2 (\xi_1 \eta_2 + \xi_2 \eta_1) + 2\lambda \xi_2 \beta_2 \eta_2) a_3 = \lambda \alpha_2 (\xi_1 \eta_2 + \xi_2 \eta_1) a_3.$

It follows that $\lambda \alpha_2 (\xi_1 \eta_2 + \xi_2 \eta_1) = 0$. If $\xi_1 = \eta_2 = \xi_2 = 1$, $\eta_1 = 0$, and then we obtain $\lambda \alpha_2 = 0$. Since $\lambda \neq 0$, $\alpha_2 = 0$. Thus, if **char**$(F) = 2$, then $\Xi(f)$ is the following matrix:

$$\begin{pmatrix} \alpha_1 & 0 & 0 \\ 0 & \beta_2 & 0 \\ \alpha_3 & \beta_3 & 0 \end{pmatrix}, \ \alpha_1, \alpha_3, \beta_2, \beta_3 \in F.$$

Now, suppose that **char**$(F) \neq 2$ and let $\mathbf{f}$ be the linear transformation of $L$ such that

$$\Xi(f) = \begin{pmatrix} \alpha_1 & 2\alpha_1 & 0 \\ \alpha_2 & \alpha_1 & 0 \\ \alpha_3 & \beta_3 & 2\alpha_1 \end{pmatrix}, \alpha_1, \alpha_2, \alpha_3, \beta_3 \in F.$$

Let $x$, $y$ be arbitrary elements of $L$, $x = \xi_1 a_1 + \xi_2 a_2 + \xi_3 a_3$, $y = \eta_1 a_1 + \eta_2 a_2 + \eta_3 a_3$ where $\xi_1, \xi_2, \xi_3, \eta_1, \eta_2, \eta_3$ are arbitrary scalars. Then

$[x, y] = (\xi_1 \eta_1 + \lambda \xi_2 \eta_2) a_3;$

$\mathbf{f}(x) = \mathbf{f}(\xi_1 a_1 + \xi_2 a_2 + \xi_3 a_3) = \xi_1 \mathbf{f}(a_1) + \xi_2 \mathbf{f}(a_2) + \xi_3 \mathbf{f}(a_3) = \xi_1 (\alpha_1 a_1 + 2\alpha_2 a_2 + \alpha_3 a_3) + \xi_2(2\alpha_1 a_1 + \alpha_1 a_2 + \beta_3 a_3) + 2\xi_3 \alpha_1 a_3 = (\xi_1 \alpha_1 + 2\xi_2 \alpha_1) a_1 + (\xi_1 \alpha_2 + \xi_2 \alpha_1) a_2 + (\xi_1 \alpha_3 + \xi_2 \beta_3 + 2\xi_3 \alpha_1) a_3;$

$\mathbf{f}(y) = (\eta_1 \alpha_1 + 2\eta_2 \alpha_1) a_1 + (\eta_1 \alpha_2 + \eta_2 \alpha_1) a_2 + (\eta_1 \alpha_3 + \eta_2 \beta_3 + 2\eta_3 \alpha_1) a_3;$

$\mathbf{f}([x, y]) = \mathbf{f}((\xi_1 \eta_1 + \lambda \xi_2 \eta_2) a_3) = (\xi_1 \eta_1 + \lambda \xi_2 \eta_2) \mathbf{f}(a_3) = 2\alpha_1 (\xi_1 \eta_1 + \lambda \xi_2 \eta_2) a_3;$

$[\mathbf{f}(x), y] + [x, \mathbf{f}(y)] = [(\xi_1 \alpha_1 + 2\xi_2 \alpha_1) a_1 + (\xi_1 \alpha_2 + \xi_2 \alpha_1) a_2 + (\xi_1 \alpha_3 + \xi_2 \beta_3 + 2\xi_3 \alpha_1) a_3, \eta_1 a_1 + \eta_2 a_2 + \eta_3 a_3] + [\xi_1 a_1 + \xi_2 a_2 + \xi_3 a_3, (\eta_1 \alpha_1 + 2\eta_2 \alpha_1) a_1 + (\eta_1 \alpha_2 + \eta_2 \alpha_1) a_2 + (\eta_1 \alpha_3 + \eta_2 \beta_3 + 2\eta_3 \alpha_1) a_3] =$

$(\xi_1\alpha_1 + 2\xi_2\alpha_1)\eta_1[a_1, a_1] + (\xi_1\alpha_2 + \xi_2\alpha_1)\eta_2[a_2, a_2] +$

$\xi_1(\eta_1\alpha_1 + 2\eta_2\alpha_1)[a_1, a_1] + \xi_2(\eta_1\alpha_2 + \eta_2\alpha_1)[a_2, a_2] =$

$(\xi_1\alpha_1 + 2\xi_2\alpha_1)\eta_1 a_3 + \lambda(\xi_1\alpha_2 + \xi_2\alpha_1)\eta_2 a_3 + \xi_1(\eta_1\alpha_1 + 2\eta_2\alpha_1)a_3 + \lambda\xi_2(\eta_1\alpha_2 + \eta_2\alpha_1)a_3 =$

$(\xi_1\eta_1\alpha_1 + 2\xi_2\alpha_1\eta_1 + \lambda\xi_1\eta_2\alpha_2 + \lambda\xi_2\eta_2\alpha_1 + \xi_1\eta_1\alpha_1 + 2\xi_1\eta_2\alpha_1 + \lambda\xi_2\eta_1\alpha_2 + \lambda\xi_2\eta_2\alpha_1)a_3 =$

$(\alpha_1 (2\xi_1\eta_1 + \lambda\xi_2\eta_2 + 2\xi_2\eta_1 + 2\xi_1\eta_2) + \lambda\alpha_2(\xi_1\eta_2 + \xi_2\eta_1))a_3.$

Hence, we obtain:

$$\alpha_1(2\xi_1\eta_1 + \lambda\xi_2\eta_2 + 2\xi_2\eta_1 + 2\xi_1\eta_2) + \lambda\alpha_2(\xi_1\eta_2 + \xi_2\eta_1) = 2\alpha_1(\xi_1\eta_1 + \lambda\xi_2\eta_2).$$

It follows that

$$\alpha_1(2\xi_2\eta_1 + 2\xi_1\eta_2 - \lambda\xi_2\eta_2) + \lambda\alpha_2(\xi_1\eta_2 + \xi_2\eta_1) = 0.$$

Put $\xi_1 = \eta_2 = \xi_2 = 1$, $\eta_1 = -1$; then we obtain $\alpha_1(-2 + 2 - \lambda) = 0$. Since $\lambda \neq 0$, $\alpha_1 = 0$. Then $\lambda\alpha_2(\xi_1\eta_2 + \xi_2\eta_1) = 0$. Put $\xi_1 = \eta_2 = \xi_2 = \eta_1 = 1$; then we obtain $2\lambda\alpha_2 = 0$. Since $\lambda \neq 0$, $\alpha_1 = 0$. Thus, if **char**$(F) \neq 2$, then $\Xi(f)$ is the following matrix:

$$\begin{pmatrix} 0 & 0 & 0 \\ 0 & 0 & 0 \\ \alpha_3 & \beta_3 & 0 \end{pmatrix}, \alpha_3, \beta_3 \in F.$$

We describe the structure of Lie algebra $\mathbb{D}$ in detail.

First, suppose that **char**$(F) = 2$. Let $\textbf{z}$ be the linear transformation of L defined by the rule $\textbf{z}(a_1) = a_3$, $\textbf{z}(a_2) = 0$, $\textbf{z}(a_3) = 0$. As we have seen above, mapping $\textbf{z}$ is a derivation of L. We note that $\Xi(\textbf{z})$ is the following matrix:

$$\Xi(\textbf{z}) = \begin{pmatrix} 0 & 0 & 0 \\ 0 & 0 & 0 \\ 1 & 0 & 0 \end{pmatrix}.$$

Denote by $\textbf{Z}$ the subspace of $\mathbb{D}$ generated by a derivation $\textbf{z}$. Clearly, $\textbf{Z}$ is a subalgebra of $\mathbb{D}$. Moreover, $\textbf{Z}$ is an ideal of $\mathbb{D}$. The following equations show this:

$$\begin{pmatrix} \alpha_1 & 0 & 0 \\ 0 & \beta_2 & 0 \\ \alpha_3 & \beta_3 & 0 \end{pmatrix}\begin{pmatrix} 0 & 0 & 0 \\ 0 & 0 & 0 \\ \sigma & 0 & 0 \end{pmatrix} - \begin{pmatrix} 0 & 0 & 0 \\ 0 & 0 & 0 \\ \sigma & 0 & 0 \end{pmatrix}\begin{pmatrix} \alpha_1 & 0 & 0 \\ 0 & \beta_2 & 0 \\ \alpha_3 & \beta_3 & 0 \end{pmatrix} =$$

$$\begin{pmatrix} 0 & 0 & 0 \\ 0 & 0 & 0 \\ 0 & 0 & 0 \end{pmatrix} - \begin{pmatrix} 0 & 0 & 0 \\ 0 & 0 & 0 \\ \sigma\alpha_1 & 0 & 0 \end{pmatrix} = \begin{pmatrix} 0 & 0 & 0 \\ 0 & 0 & 0 \\ -\sigma\alpha_1 & 0 & 0 \end{pmatrix}.$$

Let $\textbf{w}$ be the linear transformation of L defined by the rule $\textbf{w}(a_1) = 0$, $\textbf{w}(a_2) = a_3$, $\textbf{w}(a_3) = 0$. By the proof above, mapping $\textbf{w}$ is a derivation of L. We note that $\Xi(\textbf{w})$ is the following matrix:

$$\Xi(\textbf{w}) = \begin{pmatrix} 0 & 0 & 0 \\ 0 & 0 & 0 \\ 0 & 1 & 0 \end{pmatrix}.$$

Denote by $\textbf{W}$ the subspace of $\mathbb{D}$ generated by derivation $\textbf{w}$. We have:

$$\begin{pmatrix} \alpha_1 & 0 & 0 \\ 0 & \beta_2 & 0 \\ \alpha_3 & \beta_3 & 0 \end{pmatrix} \begin{pmatrix} 0 & 0 & 0 \\ 0 & 0 & 0 \\ 0 & \sigma & 0 \end{pmatrix} - \begin{pmatrix} 0 & 0 & 0 \\ 0 & 0 & 0 \\ 0 & \sigma & 0 \end{pmatrix} \begin{pmatrix} \alpha_1 & 0 & 0 \\ 0 & \beta_2 & 0 \\ \alpha_3 & \beta_3 & 0 \end{pmatrix} =$$

$$\begin{pmatrix} 0 & 0 & 0 \\ 0 & 0 & 0 \\ 0 & 0 & 0 \end{pmatrix} - \begin{pmatrix} 0 & 0 & 0 \\ 0 & 0 & 0 \\ 0 & \sigma\beta_2 & 0 \end{pmatrix} = \begin{pmatrix} 0 & 0 & 0 \\ 0 & 0 & 0 \\ 0 & -\sigma\beta_2 & 0 \end{pmatrix}.$$

These equations show that subspace **W** is an ideal of D, so that **W** + **Z** is an abelian ideal of $\mathbb{D}$.

Let **u**, **v** be the linear transformations of L defined by the rule $\mathbf{u}(a_1) = a_1$, $\mathbf{u}(a_2) = 0$, $\mathbf{u}(a_3) = 0$, $\mathbf{v}(a_1) = 0$, $\mathbf{v}(a_2) = a_2$, $\mathbf{v}(a_3) = 0$. These mappings are derivations of L, and

$$\Xi(\mathbf{u}) = \begin{pmatrix} 1 & 0 & 0 \\ 0 & 0 & 0 \\ 0 & 0 & 0 \end{pmatrix}, \; \Xi(\mathbf{v}) = \begin{pmatrix} 0 & 0 & 0 \\ 0 & 1 & 0 \\ 0 & 0 & 0 \end{pmatrix}.$$

Denote by **U** the subspace of $\mathbb{D}$ generated by derivation **u**, and by **V**, the subspace of $\mathbb{D}$ generated by derivation **v**. Clearly, **U** + **V** is an abelian subalgebra of $\mathbb{D}$. Moreover, $\mathbb{D}$ is a direct sum of abelian ideals **W** + **Z** and an abelian subalgebra. The above equations show that

$$[\mathbf{u}, \mathbf{z}] = -\mathbf{z}, [\mathbf{u}, \mathbf{w}] = -\mathbf{w}, [\mathbf{v}, \mathbf{z}] = \mathbb{O}, [\mathbf{v}, \mathbf{w}] = -\mathbf{w}.$$

Now, suppose that **char**(F) ≠ 2. Let **z**, **w** be linear transformations of L defined by the rule $\mathbf{z}(a_1) = a_3$, $\mathbf{z}(a_2) = 0$, $\mathbf{z}(a_3) = 0$, $\mathbf{w}(a_1) = 0$, $\mathbf{w}(a_2) = a_3$, $\mathbf{w}(a_3) = 0$. By what is proven above, the mappings **z**, **w** are derivations of L and

$$\Xi(\mathbf{z}) = \begin{pmatrix} 0 & 0 & 0 \\ 0 & 0 & 0 \\ 1 & 0 & 0 \end{pmatrix}, \; \Xi(\mathbf{w}) = \begin{pmatrix} 0 & 0 & 0 \\ 0 & 0 & 0 \\ 0 & 1 & 0 \end{pmatrix}.$$

Denote by **Z** the subspace of $\mathbb{D}$ generated by derivation **z**, and by **W**, the subspace of $\mathbb{D}$ generated by derivation **w**. It is not hard to prove that **Z** + **W** is an abelian subalgebra of $\mathbb{D}$, and moreover, $\mathbb{D} = \mathbf{W} + \mathbf{Z}$.

**Theorem 2.** *Let D be an algebra of derivations of Leibniz algebra $\mathbf{Lei}_5(3, F)$.*

*If **char**(F) = 2, then D is isomorphic to a Lie subalgebra of $\mathcal{M}_3(F)$ consisting of the following form's matrices:*

$$\begin{pmatrix} \alpha & 0 & 0 \\ 0 & \alpha & 0 \\ \gamma & \kappa & 0 \end{pmatrix}, \; \alpha, \gamma, \kappa \in F.$$

*Furthermore, D is the direct sum of an abelian ideal **W** + **Z** and an abelian subalgebra **U** where **Z** is the subspace of **D** generated by a derivation **z** defined by the rule **z**($a_1$) = $a_3$, **z**($a_2$) = 0, **z**($a_3$) = 0; **W** is the subspace of **D** generated by a derivation **w** defined by the rule **w**($a_1$) = 0, **w**($a_2$) = $a_3$, **w**($a_3$) = 0; **U** is the subspace of **D** generated by a derivation **u** defined by the rule **u**($a_1$) = $a_1$, **u**($a_2$) = $a_2$, **u**($a_3$) = 0. Moreover, **[u, z]** = − **z**, **[u, w]** = − **w**.*

*If **char**(F) ≠ 2, then D is isomorphic to a Lie subalgebra of **M**$_3$(F) consisting of the following form's matrices:*

$$\begin{pmatrix} 0 & 0 & 0 \\ 0 & 0 & 0 \\ \alpha & \beta & 0 \end{pmatrix}, \ \alpha, \beta \in F.$$

*Furthermore, D is the direct sum of subalgebras **W**, **Z** where **Z** is the subspace of **D** generated by a derivation **z** defined by the rule **z**($a_1$) = $a_3$, **z**($a_2$) = 0, **z**($a_3$) = 0; **W** is the subspace of **D** generated by a derivation **w** defined by the rule **w**($a_1$) = 0, **w**($a_2$) = $a_3$, **w**($a_3$) = 0; in particular, **D** is abelian.*

PROOF. Let L = **Lei**$_5$(3, F) and f ∈ **Der**(L). By Lemma 1, f(Fa$_3$) ≤ Fa$_3$, so that f($a_1$) =s $\alpha_1 a_1 + \alpha_2 a_2 + \alpha_3 a_3$, f($a_2$) = $\beta_1 a_1 + \beta_2 a_2 + \beta_3 a_3$, f($a_3$) = $\gamma a_3$. Then

f($a_3$) = f([$a_1, a_1$]) = [f($a_1$), $a_1$] + [$a_1$, f($a_1$)] =
[$\alpha_1 a_1 + \alpha_2 a_2 + \alpha_3 a_3, a_1$] + [$a_1, \alpha_1 a_1 + \alpha_2 a_2 + \alpha_3 a_3$] =
$\alpha_1[a_1, a_1] + \alpha_1[a_1, a_1] + \alpha_2[a_1, a_2] = (2\alpha_1 + \alpha_2)_3$;
f($a_3$) = f($\lambda^{-1}[a_2, a_2]$) = $\lambda^{-1}$f([$a_2, a_2$]) = $\lambda^{-1}$[f($a_2$), $a_2$]] + $\lambda^{-1}$([$a_2$, f($a_2$)] =
$\lambda^{-1}$ [$\beta_1 a_1 + \beta_2 a_2 + \beta_3 a_3, a_2$] + $\lambda^{-1}$[$a_2, \beta_1 a_1 + \beta_2 a_2 + \beta_3 a_3$] =
$\lambda^{-1}\beta_1[a_1, a_2] + \lambda^{-1}\beta_2[a_2, a_2] + \lambda^{-1}\beta_2[a_2, a_2] = \lambda^{-1}\beta_1 a_3 + 2\beta_2 a_3 = (\lambda^{-1}\beta_1 + 2\beta_2)a_3$;
f($a_3$) = f([$a_1, a_2$]) = [f($a_1$), $a_2$] + [$a_1$, f($a_2$)] =
[$\alpha_1 a_1 + \alpha_2 a_2 + \alpha_3 a_3, a_2$] + [$a_1, \beta_1 a_1 + \beta_2 a_2 + \beta_3 a_3$] =
$\alpha_1[a_1, a_2] + \alpha_2[a_2, a_2] + \beta_1[a_1, a_1] + \beta_2[a_1, a_2] = \alpha_1 a_3 + \lambda\alpha_2 a_3 + \beta_1 a_3 + \beta_2 a_3 =$
$(\alpha_1 + \lambda\alpha_2 + \beta_1 + \beta_2)a_3$.

Thus, we obtain the following equalities:

$2\alpha_1 + \alpha_2 = \lambda^{-1}\beta_1 + 2\beta_2 = \alpha_1 + \lambda\alpha_2 + \beta_1 + \beta_2$.

If **char**(F) = 2, then we obtain $\alpha_2 = \lambda^{-1}\beta_1 = \alpha_1 + \lambda\alpha_2 + \beta_1 + \beta_2$. It follows that $\beta_1 = \lambda\alpha_2$ and $\alpha_2 = \alpha_1 + \lambda\alpha_2 + \lambda\alpha_2 + \beta_2 = \alpha_1 + \beta_2$, so that and $\beta_2 = \alpha_2 - \alpha_1 = \alpha_2 + \alpha_1$. Thus, Ξ(f) is the following matrix:

$$\begin{pmatrix} \alpha_1 & \lambda\alpha_2 & 0 \\ \alpha_2 & \alpha_2 + \alpha_1 & 0 \\ \alpha_3 & \beta_3 & 0 \end{pmatrix}, \ \alpha_1, \alpha_2, \alpha_3, \beta_3 \in F.$$

If **char**(F) ≠ 2, then Ξ(f) is the following matrix:

$$\begin{pmatrix} \alpha_1 & 2\alpha_1 & 0 \\ \alpha_2 & \alpha_1 & 0 \end{pmatrix}$$

$\begin{pmatrix} \alpha_3 & \beta_3 & 2\alpha_1 \end{pmatrix}$, $\alpha_1, \alpha_2, \alpha_3, \beta_3 \in F$.

Conversely, suppose that $\mathbf{char}(F) = 2$ and let $\mathbf{f}$ be a linear transformation of L such that

$$\Xi(f) = \begin{pmatrix} \alpha_1 & \lambda\alpha_2 & 0 \\ \alpha_2 & \alpha_2 + \alpha_1 & 0 \\ \alpha_3 & \beta_3 & 0 \end{pmatrix}, \; \alpha_1, \alpha_2, \alpha_3, \beta_3 \in F.$$

Let x, y be the arbitrary elements of L, $x = \xi_1 a_1 + \xi_2 a_2 + \xi_3 a_3$, $y = \eta_1 a_1 + \eta_2 a_2 + \eta_3 a_3$ where $\xi_1, \xi_2, \xi_3, \eta_1, \eta_2, \eta_3$ are arbitrary scalars. Then

$[x, y] = [\xi_1 a_1 + \xi_2 a_2 + \xi_3 a_3, \eta_1 a_1 + \eta_2 a_2 + \eta_3 a_3] =$
$[\xi_1 a_1 + \xi_2 a_2 + \xi_3 a_3, \eta_1 a_1 + \eta_2 a_2 + \eta_3 a_3] = \xi_1 \eta_1 [a_1, a_1] + \xi_1 \eta_2 [a_1, a_2] + \xi_2 \eta_2 [a_2, a_2] =$
$\xi_1 \eta_1 a_3 + \xi_1 \eta_2 a_3 + \lambda \xi_2 \eta_2 a_3 = (\xi_1 \eta_1 + \xi_1 \eta_2 + \lambda \xi_2 \eta_2) a_3;$
$\mathbf{f}(x) = \mathbf{f}(\xi_1 a_1 + \xi_2 a_2 + \xi_3 a_3) = \xi_1 \mathbf{f}(a_1) + \xi_2 \mathbf{f}(a_2) + \xi_3 \mathbf{f}(a_3) =$
$\xi_1(\alpha_1 a_1 + \alpha_2 a_2 + \alpha_3 a_3) + \xi_2(\lambda\alpha_2 a_1 + (\alpha_2 + \alpha_1)a_2 + \beta_3 a_3) =$
$\xi_1 \alpha_1 a_1 + \xi_1 \alpha_2 a_2 + \xi_1 \alpha_3 a_3 + \lambda \xi_2 \alpha_2 a_1 + \xi_2(\alpha_2 + \alpha_1)a_2 + \xi_2 \beta_3 a_3 =$
$(\xi_1 \alpha_1 a + \lambda \xi_2 \alpha_2)a_1 + (\xi_1 \alpha_2 + \xi_2 \alpha_2 + \xi_2 \alpha_1)a_2 + (\xi_1 \alpha_3 + \xi_2 \beta_3)a_3;$
$\mathbf{f}(y) = (\eta_1 \alpha_1 a + \lambda \eta_2 \alpha_2)a_1 + (\eta_1 \alpha_2 + \eta_2 \alpha_2 + \eta_2 \alpha_1)a_2 + (\eta_1 \alpha_3 + \eta_2 \beta_3)a_3;$
$\mathbf{f}([x, y]) = \mathbf{f}((\xi_1 \eta_1 + \xi_1 \eta_2 + \lambda \xi_2 \eta_2)a_3) = (\xi_1 \eta_1 + \xi_1 \eta_2 + \lambda \xi_2 \eta_2)\mathbf{f}(a_3) = 0,$
$[\mathbf{f}(x), y] + [x, \mathbf{f}(y)] =$
$[(\xi_1 \alpha_1 a + \lambda \xi_2 \alpha_2)a_1 + (\xi_1 \alpha_2 + \xi_2 \alpha_2 + \xi_2 \alpha_1)a_2 + (\xi_1 \alpha_3 + \xi_2 \beta_3)a_3, \eta_1 a_1 + \eta_2 a_2 + \eta_3 a_3] +$
$[\xi_1 a_1 + \xi_2 a_2 + \xi_3 a_3, (\eta_1 \alpha_1 a + \lambda \eta_2 \alpha_2)a_1 + (\eta_1 \alpha_2 + \eta_2 \alpha_2 + \eta_2 \alpha_1)a_2 + (\eta_1 \alpha_3 + \eta_2 \beta_3)a_3] =$
$(\xi_1 \alpha_1 + \lambda \xi_2 \alpha_2)\eta_1 [a_1, a_1] + (\xi_1 \alpha_1 a + \lambda \xi_2 \alpha_2)\eta_2 [a_1, a_2] + (\xi_1 \alpha_2 + \xi_2 \alpha_2 + \xi_2 \alpha_1)\eta_2 [a_2, a_2] +$
$\xi_1(\eta_1 \alpha_1 a + \lambda \eta_2 \alpha_2)[a_1, a_1] + \xi_1(\eta_1 \alpha_2 + \eta_2 \alpha_2 + \eta_2 \alpha_1)[a_1, a_2] + \xi_2(\eta_1 \alpha_2 + \eta_2 \alpha_2 + \eta_2 \alpha_1)[a_2, a_2]$
$= (\xi_1 \alpha_1 + \lambda \xi_2 \alpha_2)\eta_1 a_3 + (\xi_1 \alpha_1 + \lambda \xi_2 \alpha_2)\eta_2 a_3 + \lambda(\xi_1 \alpha_2 + \xi_2 \alpha_2 + \xi_2 \alpha_1)\eta_2 a_3 +$
$\xi_1(\eta_1 \alpha_1 + \lambda \eta_2 \alpha_2)a_3 + \xi_1(\eta_1 \alpha_2 + \eta_2 \alpha_2 + \eta_2 \alpha_1)a_3 + \lambda \xi_2(\eta_1 \alpha_2 + \eta_2 \alpha_2 + \eta_2 \alpha_1)a_3 =$
$(\xi_1 \alpha_1 \eta_1 + \lambda \xi_2 \alpha_2 \eta_1 + \xi_1 \alpha_1 \eta_2 + \lambda \xi_2 \alpha_2 \eta_2 + \lambda \xi_1 \alpha_2 \eta_2 + \lambda \xi_2 \alpha_2 \eta_2 + \lambda \xi_2 \alpha_1 \eta_2 +$
$\xi_1 \eta_1 \alpha_1 + \lambda \xi_1 \eta_2 \alpha_2 + \xi_1 \eta_1 \alpha_2 + \xi_1 \eta_2 \alpha_2 + \xi_1 \eta_2 \alpha_1 + \lambda \xi_2 \eta_1 \alpha_2 + \lambda \xi_2 \eta_2 \alpha_2 + \lambda \xi_2 \eta_2 \alpha_1)a_3 =$
$(\xi_1 \eta_1 \alpha_2 + \xi_1 \eta_2 \alpha_2 + \lambda \xi_2 \eta_2 \alpha_2)a_3 = \alpha_2(\xi_1 \eta_1 + \xi_1 \eta_2 + \lambda \xi_2 \eta_2)a_3.$

Thus, we obtain $\alpha_2(\xi_1 \eta_1 + \xi_1 \eta_2 + \lambda \xi_2 \eta_2) = 0$. If $\xi_1 = \eta_1 = \xi_2 = 1$, $\eta_2 = 0$, then we obtain $\alpha_2 = 0$. Hence, if $\mathbf{char}(F) = 2$, then $\Xi(f)$ is the following matrix:

$$\Xi(f) = \begin{pmatrix} \alpha_1 & 0 & 0 \\ 0 & \alpha_1 & 0 \\ \alpha_3 & \beta_3 & 0 \end{pmatrix}, \; \alpha_1, \alpha_3, \beta_3 \in F.$$

Suppose that $\mathbf{char}(F) \neq 2$ and let $\mathbf{f}$ be a linear transformation of L such that

$$\Xi(f) = \begin{pmatrix} \alpha_1 & 2\alpha_1 & 0 \\ \alpha_2 & \alpha_1 & 0 \\ \alpha_3 & \beta_3 & 2\alpha_1 \end{pmatrix}, \; \alpha_1, \alpha_2, \alpha_3, \beta_3 \in F.$$

Let x, y be the arbitrary elements of L, $x = \xi_1 a_1 + \xi_2 a_2 + \xi_3 a_3$, $y = \eta_1 a_1 + \eta_2 a_2 + \eta_3 a_3$ where $\xi_1, \xi_2, \xi_3, \eta_1, \eta_2, \eta_3$ are arbitrary scalars. Then

$[x, y] = [\xi_1 a_1 + \xi_2 a_2 + \xi_3 a_3, \eta_1 a_1 + \eta_2 a_2 + \eta_3 a_3] = (\xi_1 \eta_1 + \xi_1 \eta_2 + \lambda \xi_2 \eta_2) a_3;$

$\mathbf{f}(x) = \mathbf{f}(\xi_1 a_1 + \xi_2 a_2 + \xi_3 a_3) = \xi_1 \mathbf{f}(a_1) + \xi_2 \mathbf{f}(a_2) + \xi_3 \mathbf{f}(a_3) =$

$\xi_1(\alpha_1 a_1 + \alpha_2 a_2 + \alpha_3 a_3) + \xi_2(2\alpha_1 a_1 + \alpha_1 a_2 + \beta_3 a_3) + 2\xi_3 \alpha_1 a_3 =$

$\xi_1 \alpha_1 a_1 + \xi_1 \alpha_2 a_2 + \xi_1 \alpha_3 a_3 + 2\xi_2 \alpha_1 a_1 + \xi_2 \alpha_1 a_2 + \xi_2 \beta_3 a_3 + 2\xi_3 \alpha_1 a_3 =$

$(\xi_1 \alpha_1 + 2\xi_2 \alpha_1) a_1 + (\xi_1 \alpha_2 + \xi_2 \alpha_1) a_2 + (\xi_1 \alpha_3 + \xi_2 \beta_3 + 2\xi_3 \alpha_1) a_3;$

$\mathbf{f}(y) = (\eta_1 \alpha_1 + 2\eta_2 \alpha_1) a_1 + (\eta_1 \alpha_2 + \eta_2 \alpha_1) a_2 + (\eta_1 \alpha_3 + \eta_2 \beta_3 + 2\eta_3 \alpha_1) a_3;$

$\mathbf{f}([x, y]) = \mathbf{f}((\xi_1 \eta_1 + \xi_1 \eta_2 + \lambda \xi_2 \eta_2) a_3) = 2\alpha_1 (\xi_1 \eta_1 + \xi_1 \eta_2 + \lambda \xi_2 \eta_2);$

$[\mathbf{f}(x), y] + [x, \mathbf{f}(y)] =$

$[(\xi_1 \alpha_1 + 2\xi_2 \alpha_1) a_1 + (\xi_1 \alpha_2 + \xi_2 \alpha_1) a_2 + (\xi_1 \alpha_3 + \xi_2 \beta_3 + 2\xi_3 \alpha_1) a_3, \eta_1 a_1 + \eta_2 a_2 + \eta_3 a_3] +$

$[\xi_1 a_1 + \xi_2 a_2 + \xi_3 a_3, (\eta_1 \alpha_1 + 2\eta_2 \alpha_1) a_1 + (\eta_1 \alpha_2 + \eta_2 \alpha_1) a_2 + (\eta_1 \alpha_3 + \eta_2 \beta_3 + 2\eta_3 \alpha_1) a_3] =$

$(\xi_1 \alpha_1 + 2\xi_2 \alpha_1) \eta_1 [a_1, a_1] + (\xi_1 \alpha_1 + 2\xi_2 \alpha_1) \eta_2 [a_1, a_2] + (\xi_1 \alpha_2 + \xi_2 \alpha_1) \eta_2 [a_2, a_2] +$

$\xi_1 (\eta_1 \alpha_1 + 2\eta_2 \alpha_1) [a_1, a_1] + \xi_1 (\eta_1 \alpha_2 + \eta_2 \alpha_1) [a_1, a_2] + \xi_2 (\eta_1 \alpha_2 + \eta_2 \alpha_1) [a_2, a_2] =$

$(\xi_1 \alpha_1 + 2\xi_2 \alpha_1) \eta_1 a_3 + (\xi_1 \alpha_1 + 2\xi_2 \alpha_1) \eta_2 a_3 + \lambda (\xi_1 \alpha_2 + \xi_2 \alpha_1) \eta_2 a_3 + \xi_1 (\eta_1 \alpha_1 + 2\eta_2 \alpha_1) a_3 +$

$\xi_1 (\eta_1 \alpha_2 + \eta_2 \alpha_1) a_3 + \lambda \xi_2 (\eta_1 \alpha_2 + \eta_2 \alpha_1) a_3 =$

$(\xi_1 \alpha_1 \eta_1 + 2\xi_2 \alpha_1 \eta_1 + \xi_1 \alpha_1 \eta_2 + 2\xi_2 \alpha_1 \eta_2 + \lambda \xi_1 \alpha_2 \eta_2 + \lambda \xi_2 \alpha_1 \eta_2 + \xi_1 \eta_1 \alpha_1 + 2\xi_1 \eta_2 \alpha_1 +$

$\xi_1 \eta_1 \alpha_2 + \xi_1 \eta_2 \alpha_1 + \lambda \xi_2 \eta_1 \alpha_2 + \lambda \xi_2 \eta_2 \alpha_1) a_3 =$

$(2\alpha_1 (\xi_1 \eta_1 + \xi_2 \eta_1 + \xi_2 \eta_2 + \xi_1 \eta_2 + \lambda \xi_2 \eta_2) + \alpha_2 (\lambda \xi_1 \eta_2 + \xi_1 \eta_1 + \lambda \xi_2 \eta_1) a_3.$

It follows that

$2\alpha_1 (\xi_1 \eta_1 + \xi_2 \eta_1 + \xi_2 \eta_2 + \xi_1 \eta_2 + \lambda \xi_2 \eta_2) + \alpha_2 (\lambda \xi_1 \eta_2 + \xi_1 \eta_1 + \lambda \xi_2 \eta_1) =$
$2\alpha_1 (\xi_1 \eta_1 + \xi_1 \eta_2 + \lambda \xi_2 \eta_2).$

Then we obtain:

$2\alpha_1 \xi_2 (\eta_1 + \eta_2) + \alpha_2 (\lambda \xi_1 \eta_2 + \xi_1 \eta_1 + \lambda \xi_2 \eta_1) = 0.$

If $\xi_1 = \eta_1 = 0$, $\xi_2 = \eta_2 = 1$, then we will have $2\alpha_1 = 0$, and hence, $\alpha_1 = 0$. In this case, $\alpha_2 (\lambda \xi_1 \eta_2 + \xi_1 \eta_1 + \lambda \xi_2 \eta_1) = 0$. If $\xi_1 = 0$, then $\alpha_2 \lambda \xi_2 \eta_1 = 0$, and it follows that $\alpha_2 = 0$. Hence, if $\mathbf{char}(F) \neq 2$, then $\Xi(f)$ is the following matrix:

$$\Xi(f) = \begin{pmatrix} 0 & 0 & 0 \\ 0 & 0 & 0 \\ \alpha_3 & \beta_3 & 0 \end{pmatrix}, \quad \alpha_3, \beta_3 \in F.$$

We describe the structure of Lie algebra $\mathbb{D}$ in detail.

Suppose first that $\mathbf{char}(F) = 2$. Let $\boldsymbol{z}$ be a linear transformation of L defined by the rule $\boldsymbol{z}(a_1) = a_3$, $\boldsymbol{z}(a_2) = 0$, $\boldsymbol{z}(a_3) = 0$. As we have seen above, mapping $\boldsymbol{z}$ is a derivation of L. We note that $\Xi(\boldsymbol{z})$ is the following matrix:

$$\Xi(\boldsymbol{z}) = \begin{pmatrix} 0 & 0 & 0 \\ 0 & 0 & 0 \\ 1 & 0 & 0 \end{pmatrix}.$$

Denote by $\boldsymbol{Z}$ the subspace of $\mathbb{D}$, generated by a(the?) derivation $\boldsymbol{z}$. Clearly, $\boldsymbol{Z}$ is a subalgebra of $\mathbb{D}$. Moreover, $\boldsymbol{Z}$ is an ideal of $\mathbb{D}$. The following equations show it.

$$\begin{pmatrix} \alpha_1 & 0 & 0 \\ 0 & \alpha_1 & 0 \\ \alpha_3 & \beta_3 & 0 \end{pmatrix} \begin{pmatrix} 0 & 0 & 0 \\ 0 & 0 & 0 \\ \sigma & 0 & 0 \end{pmatrix} - \begin{pmatrix} 0 & 0 & 0 \\ 0 & 0 & 0 \\ \sigma & 0 & 0 \end{pmatrix} \begin{pmatrix} \alpha_1 & 0 & 0 \\ 0 & \alpha_1 & 0 \\ \alpha_3 & \beta_3 & 0 \end{pmatrix} =$$

$$\begin{pmatrix} 0 & 0 & 0 \\ 0 & 0 & 0 \\ 0 & 0 & 0 \end{pmatrix} - \begin{pmatrix} 0 & 0 & 0 \\ 0 & 0 & 0 \\ \sigma\alpha_1 & 0 & 0 \end{pmatrix} = \begin{pmatrix} 0 & 0 & 0 \\ 0 & 0 & 0 \\ -\sigma\alpha_1 & 0 & 0 \end{pmatrix}.$$

Let $\mathbf{w}$ be a linear transformation of L defined by the rule $\mathbf{w}(a_1) = 0$, $\mathbf{w}(a_2) = a_3$, $\mathbf{w}(a_3) = 0$. By what is proven above, mapping $\mathbf{w}$ is a derivation of L. We note that $\Xi(\mathbf{w})$ is the following matrix:

$$\Xi(\mathbf{w}) = \begin{pmatrix} 0 & 0 & 0 \\ 0 & 0 & 0 \\ 0 & 1 & 0 \end{pmatrix}.$$

Denote by $\mathbf{W}$ the subspace of $\mathbb{D}$ generated by derivation $\mathbf{w}$. We have:

$$\begin{pmatrix} \alpha_1 & 0 & 0 \\ 0 & \alpha_1 & 0 \\ \alpha_3 & \beta_3 & 0 \end{pmatrix} \begin{pmatrix} 0 & 0 & 0 \\ 0 & 0 & 0 \\ 0 & \sigma & 0 \end{pmatrix} - \begin{pmatrix} 0 & 0 & 0 \\ 0 & 0 & 0 \\ 0 & \sigma & 0 \end{pmatrix} \begin{pmatrix} \alpha_1 & 0 & 0 \\ 0 & \alpha_1 & 0 \\ \alpha_3 & \beta_3 & 0 \end{pmatrix} =$$

$$\begin{pmatrix} 0 & 0 & 0 \\ 0 & 0 & 0 \\ 0 & 0 & 0 \end{pmatrix} - \begin{pmatrix} 0 & 0 & 0 \\ 0 & 0 & 0 \\ 0 & \sigma\beta_2 & 0 \end{pmatrix} = \begin{pmatrix} 0 & 0 & 0 \\ 0 & 0 & 0 \\ 0 & -\sigma\alpha_1 & 0 \end{pmatrix}.$$

These equations show that subspace $\mathbf{W}$ is an ideal of D, so that $\mathbf{W} + \mathbf{Z}$ is an abelian ideal of $\mathbb{D}$.

Let $\mathbf{u}$ be a linear transformation of L defined by the rule $\mathbf{u}(a_1) = a_1$, $\mathbf{u}(a_2) = a_2$, $\mathbf{u}(a_3) = 0$. This mapping is the derivations of L, and

$$\Xi(\mathbf{u}) = \begin{pmatrix} 1 & 0 & 0 \\ 0 & 1 & 0 \\ 0 & 0 & 0 \end{pmatrix}.$$

Denote by $\mathbb{U}$ a subspace of $\mathbb{D}$ generated by derivation $\mathbf{u}$. Clearly, $\mathbb{U}$ is an abelian subalgebra of $\mathbb{D}$. Moreover, $\mathbb{D}$ is a direct sum of abelian ideals $\mathbf{W} + \mathbf{Z}$ and this abelian subalgebra $\mathbb{U}$. The equations above show that

$$[\mathbf{u}, \mathbf{z}] = -\mathbf{z}, \ [\mathbf{u}, \mathbf{w}] = = -\mathbf{w}.$$

Now, suppose that $\mathbf{char}(F) \neq 2$. Let $\mathbf{z}, \mathbf{w}$ be the linear transformations of L defined by the rule $\mathbf{z}(a_1) = a_3$, $\mathbf{z}(a_2) = 0$, $\mathbf{z}(a_3) = 0$, $\mathbf{w}(a_1) = 0$, $\mathbf{w}(a_2) = a_3$, $\mathbf{w}(a_3) = 0$. By what is proven above, the mappings $\mathbf{z}, \mathbf{w}$ are derivations of L and

$$\Xi(\mathbf{z}) = \begin{pmatrix} 0 & 0 & 0 \\ 0 & 0 & 0 \\ 1 & 0 & 0 \end{pmatrix}, \ \Xi(\mathbf{w}) = \begin{pmatrix} 0 & 0 & 0 \\ 0 & 0 & 0 \\ 0 & 1 & 0 \end{pmatrix}.$$

Denote by $\mathbb{Z}$ the subspace of $\mathbb{D}$ generated by derivation $z$ and by the subspace $\mathbb{W}$ of $\mathbb{D}$ generated by derivation $w$. It is not hard to prove that $\mathbb{Z} + \mathbb{W}$ is an abelian subalgebra of $\mathbb{D}$. Moreover, $\mathbb{D} = \mathbb{W} + \mathbb{Z}$.